\RequirePackage{fix-cm}
\documentclass{svjour3}                     
\smartqed  
\usepackage{graphicx}
\usepackage{color}
\usepackage[utf8]{inputenc}
\usepackage{amsmath}
\usepackage{amssymb,latexsym}
\usepackage{epstopdf}
\usepackage{booktabs}
\usepackage{amsfonts}
\usepackage{tikz}
\usepackage{array}
\usepackage{mathtools}
\usepackage{eucal}
\usepackage[T1]{fontenc}
\usepackage{float}
\usepackage{subfig}
\usepackage{multirow}
\usepackage{adjustbox}
\numberwithin{equation}{section}
\usepackage{tabularx,ragged2e,booktabs,caption}
\usetikzlibrary{shapes.geometric, arrows}
\usepackage{hyperref}
\usepackage{mathptmx}    

\newcommand{\B}[1]{\mbox{\boldmath${#1}$\unboldmath}}

\newcommand{\Reals}{\mathbb{R}}

\newcommand{\overbar}[1]{\mkern 1.5mu\overline{\mkern-1.5mu#1\mkern-1.5mu}\mkern 1.5mu}

\newcommand{\norm}[1]{\left\lVert#1\right\rVert}

\begin{document}

\title{Accounting for model error in Tempered Ensemble Transform Particle Filter 
and its application to non-additive model error}

\titlerunning{Model error in TETPF}        

\author{Svetlana Dubinkina        \and
        Sangeetika Ruchi
}


\institute{ S. Dubinkina \at
               CWI, Science park 123, 1098XG, Amsterdam, the Netherlands\\
              \email{s.dubinkina@cwi.nl}           
               \and
               S. Ruchi \at
              CWI, Science park 123, 1098XG, Amsterdam, the Netherlands
}

\date{Received: date / Accepted: date}

\maketitle
\begin{abstract}
In this paper, we trivially extend Tempered (Localized) Ensemble Transform Particle Filter---T(L)ETPF---to account for model error. 
We examine T(L)ETPF performance for non-additive model error
in a low-dimensional and a high-dimensional test problem. The former one is a nonlinear toy model, where 
uncertain parameters are non-Gaussian distributed but model error is Gaussian distributed.
The latter one is a steady-state single-phase Darcy flow model,
where uncertain parameters are Gaussian distributed but model error is non-Gaussian distributed.
The source of model error in the Darcy flow problem is uncertain boundary conditions.
We comapare T(L)ETPF to a Regularized (Localized) Ensemble Kalman 
Filter---R(L)EnKF. We show that T(L)ETPF outperforms R(L)EnKF for both the low-dimensional and the high-dimensional problem.
This holds even when ensemble size of TLETPF is 100 while ensemble size of R(L)EnKF is greater than 6000. 
As a side note, we show that TLETPF takes less iterations than TETPF, which decreases computational costs;
while RLEnKF takes more iterations than REnKF, which incerases computational costs.
This is due to an influence of localization on a tempering and a regularizing parameter.

\end{abstract}

\section{Introduction}

Ensemble-based data assimilation deals with 
estimation of uncertain parameters and states of a model constrained by available observations using an ensemble. 
It is widely employed in many fields, for example meteorology~\cite{Hoetal96} and reservoir engineering~\cite{EmRe13}.
While in meteorology one is interested in estimation
of uncertain initial conditions of a high-dimensional chaotic system, in reservoir engineering---of estimating high-dimensional
uncertain parameters, of permeability for example, of a deterministic non-chaotic system. 

However, uncertainty is not only in initial conditions or random coefficients of a PDE
but also in a model itself, hence model error. 
In inverse problems typical sources of model error are model reduction, when a complex model is replaced by a simple
one, and incorrect parametrization.
It has been acknowledged, i.e~\cite{KaSo04}, that accounting
for model error in data assimilation greatly improves parameter-state estimation. 
Recent advances in accounting for model error in ensemble-based data assimilation are 
extension of iterative ensemble Kalman filter to include additive model error~\cite{Saetal18},
Gaussian model error update for Bayesian inversion~\cite{Caetal18}, and adding model error in the 
randomized maximum likelihood though to correctly sample the posterior without marginalization~\cite{Oliver17}.

However, most of these works have considered either additive model error or Gaussian model error,
with the sole exception of~\cite{Oliver17} where Gaussian anamorphosis was used.
However, for high-dimensional problems finding a transformation to multivariate Gaussian probability is computationally challenging.
Though the assumption of additive Gaussian model error
simplifies an optimization problem, model error is not limited to being additive nor Gaussian.
Therefore, it is essential for a data assimilation method to account for model error in a most general way.
A straightforward example of such a data assimilation method is Markov Chain Monte Carlo (MCMC). 
However, for high-dimensional systems it is impractical.

An alternative to MCMC is particle filtering~\cite{Doucet01}. Particle filtering is based on proposing
an ensemble from a prior that is not necessary close to the target posterior and to correct 
for this mismatch my computing importance weights. The importance weights are defined as a function of 
ensemble estimations and available observations. The ensemble is then resampled according to 
the estimated posterior. 
Particle filtering in its original form worked only 
for low-dimensional problems. 
However, due to recent advances of employing localization~\cite{PeMi16,Poterjoy16} it has been proven to strive for high-dimensional
problems as well. 

There are different approaches to resampling in particle filtering, but most of them are stochastic.
An ensemble transform particle filter~\cite{ReCo15} employs {\it deterministic} resampling, 
which reduces sampling error and thus needs a smaller ensemble than a typical particle filter. 
It also has a localized version.  
In~\cite{RuDu18}, we have employed the method to an inverse problem of uncertain permeability.
We have shown that though localization makes the ensemble transform particle filter deteriorate a posterior estimation of
the leading modes, it makes the method applicable to high-dimensional problems.
In~\cite{RuDuIg18}, instead of localization we have implemented tempering to the ensemble transform particle filter (TETPF).
We have shown that iterations based on temperatures~\cite{Neal01,KaBeJa14} handle notably strongly nonlinear cases
and that TETPF is able to predict multimodal distributions for high-dimensional problems.

In the current work, we account for model error in T(L)ETPF (with and without localization).
We should note that it has been already accounted for model error in the ensemble transform particle filter in~\cite{ReCo15}.
Thus it is rather trivial extension of the ensemble transform particle filter. However, our goal here is
to investigate T(L)ETPF performance in case of non-additive non-Gaussian model error 
and to compare it to an ensemble Kalman filter. 

Ensemble Kalman filter assumes Gaussian probabilities. It solves an optimization problem for mean
and approximates variance with an ensemble. It has been shown, also in presence of model error in~\cite{Evensen18,Du18}, that 
ensemble Kalman filter is able to estimate skewed probabilities.
It, however, fails to estimate multimodal probabilities. This was shown for example in~\cite{RuDuIg18}, though without model error. 
There a regularized ensemble Kalman filter (REnKF)~\cite{Iglesias15} poorly estimated geometrical parameters of a model of a Darcy flow over a channelized domain. 

Here we consider two test problems, both of them with model error, 
and compare T(L)ETPF to R(L)EnKF (with and without localization).
First test problem is a 2D scalar model with multiplicative model error.
Uncertain parameters are multimodal distributed and model error is Gaussian distributed.
Second test problem is a steady-state single-phase Darcy flow model with model error in
uncertain boundary conditions. 
This groundwater model, thought without model error, was used
first used as benchmark for inverse modeling in~\cite{CaNe86}.
It has been also used as a test model for the identification of parameters with iterative regularization
methods in~\cite{Hanke97} and with an ensemble Kalman approach in~\cite{Iglesias16}, for example.

Model error we introduce in the groundwater model in boundary conditions makes model error non-additive.
Moreover, we make it non-Gaussian by construction.
Uncertain log permeability, however, is defined as a Gaussian process. It is well known
that ensemble Kalman filter (iterative or regualrized) gives fine estimations of Gaussian probabilities even for nonlinear models. 
Our goal, however, is to investigate whether an estimation of R(L)EnKF is sensitive to uncertain boundary conditions.
Thus whether R(L)EnKF gives worse estimation than T(L)ETPF.

%
%

\subsection{Bayesian Inference}
Both T(L)ETPF and R(L)EnKF are based on Bayesian inference.
Assume $\B{u}\in \tilde{\mathcal{U}}$ is a random variable and $\B{y}_{\rm obs}\in \mathcal{Y}$ is an observation
of $\B{u}$. Then according to the Bayes's formula
\[
 \pi(\B{u}|\B{y}_{\rm obs}) \propto \pi(\B{y}_{\rm obs}|\B{u})\pi(\B{u})
\]
up to a constant of normalization.
For any smooth function $f\colon \tilde{\mathcal{U}} \to \tilde{\mathcal{U}}$, its expectation is defined as
\[
 \overbar{f(\B{u})} = \int d\B{u} f(\B{u})\pi(\B{u}|\B{y}_{\rm obs}).
\]
Assume $\B{q}\in \tilde{\mathcal{Q}}$ is a random variable independent of $\B{u}$
and the marginal distribution is 
\begin{equation}\label{eq:marg}
 \pi(\B{u}|\B{y}_{\rm obs}) = \int d\B{q} \pi(\B{u},\B{q}|\B{y}_{\rm obs}),
\end{equation}
then the expectation of $f(\B{u})$ is
\begin{equation}\label{eq:exp}
 \overbar{f(\B{u})} = \int d\B{u} f(\B{u})\int d\B{q} \pi(\B{u},\B{q}|\B{y}_{\rm obs}),
\end{equation}
where 
\begin{equation}\label{eq:Bayes}
 \pi(\B{u},\B{q}|\B{y}_{\rm obs}) \propto \pi(\B{y}_{\rm obs}|\B{u},\B{q})\pi(\B{u})\pi(\B{q})
\end{equation}
according to the Bayes formula.

The random variable $\B{u}$ denotes uncertain parameters and the random variable $q$ denotes model error. Let us first consider a case without model error.
We denote by $G\colon\tilde{\mathcal{U}} \to \mathcal{Y}$ the nonlinear forward operator
that arises from a model under consideration. In other words,
$G$ maps the space $\tilde{\mathcal{U}}$ of uncertain quantities (states or parameters)
to the observation space $\mathcal{Y}$ defined in terms of observable quantities, which are related to the solution of the model as 
\[
 y = G(u). 
\]
Then the conditional probability density function is
\[
 \pi(\B{y}_{\rm obs}|\B{u}) = \pi(\B{y}_{\rm obs}|y)\delta(y-G(u)),
\]
where $\delta$ is the Dirac delta function. \par

When model error is present, the nonlinear forward operator is $g\colon \tilde{\mathcal{U}}\times \tilde{\mathcal{Q}} \to \mathcal{Y}$ and
related to the solution of the model as  
\[
 y = g(u,q). 
\]
Assuming the forward mapping $g$ is an erroneous approximation of the true mapping $G$, 
the joint probability density function is expressed as
\[
 \pi(\B{y}_{\rm obs}|u,q) = \pi(\B{y}_{\rm obs}|y)\pi(g(u,q)|G(u)).
\]
The disadvantage of this approach is that one has to make an assumption about the transition density $\pi(g(u,q)|G(u))$. 
For example, in~\cite{Caetal18} the transition density was assumed to be Gaussian with unknown mean and covariance.

We, instead, express the joint probability density function as
\begin{equation}\label{eq:joint}
 \pi(\B{y}_{\rm obs}|\B{u},\B{q}) = \pi(\B{y}_{\rm obs}|y)\delta(y-g(u,q)),
\end{equation}
where the transition density is the Dirac delta function. 
This approach was taken in~\cite{Evensen18} for iterative ensemble smoothers to account for model error. The model error was estimated 
for both linear and nonlinear toy problems and additive model errors. 
We, on the contrary, do not 
restrict ourselves to additive model errors.

\section{Tempered Ensemble Transform Particle Filter}\label{sec:TETPF}
The goal of the Bayesian approach is to compute the posterior given by Eq.~\eqref{eq:Bayes}--\eqref{eq:joint}.
Sequential Monte Carlo (SMC) is an approximation of the Bayesian posterior.
An SMC method creates a finite sample from a prior, that is easy to sample from, and
corrects for the differences between the prior and the posterior by computing so-called importance weights. 
Finally, a resampling is performed according to those weights in order to create a new sample.
\subsection{Importance weights}

We consider discrete random variables and define $\mathcal{U} = \{\vec{u}_1,\dots, \vec{u}_M\}\subset\tilde{\mathcal{U}}$, $\vec{u}_i\in\Reals^n$.
The model has an unknown quantity $\vec{u}^{\rm true}\in \Reals^n$ that we wish to estimate from noisy observations $\vec{y}_{\rm obs}\in \Reals^\kappa$, 
where $\kappa<n$, 
\[
 \vec{y}_{\rm obs} := G(\vec{u}^{\rm true}) + \eta, 
\]
where $\eta \sim \mathcal{N}(\vec{0},\vec{R})$ with $\vec{R}$ being a known covariance matrix of the observation noise. 
Then the conditional probability density function is 
\[
 \pi(\vec{y}_{\rm obs}|\vec{y}) \propto \exp\left[-\frac{1}{2}(\vec{y}-\vec{y}_{\rm obs})^T \vec{R}^{-1}(\vec{y}-\vec{y}_{\rm obs})\right].
\]
Define $\mathcal{Q} = \{\vec{q}_1,\dots, \vec{q}_M\}\subset\tilde{\mathcal{Q}}$, $\vec{q}_i\in\Reals^m$, then 
a discrete approximation to Eq.~\eqref{eq:joint} is
\[
\pi(\vec{y}_{\rm obs}|\vec{u},\vec{q}) \propto \exp\left[-\frac{1}{2}\left(g(\vec{u},\vec{q})-\vec{y}_{\rm obs}\right)^\prime \vec{R}^{-1} \left(g(\vec{u},\vec{q})-\vec{y}_{\rm obs}\right)\right],
\]
where $\prime$ denotes the transpose.
We assume the priors $\pi(\vec{u})$ and $\pi(\vec{q})$ are uniform, then denoting $\vec{v}= [\vec{u}\ \vec{q}]^\prime$ the expectation of a function $f$ of  $\vec{v}$  is
\[
 \overbar{f(\vec{v})} \approx \sum_{i=1}^{M}f(\vec{v}_i)w_i.
 \]
Here the importance weights are
\begin{equation}\label{eq:weights}
 w_i = \frac{h(\vec{v}_i)}{\sum_{j=1}^M h(\vec{v}_j)},\quad \mbox{where} \quad
h(\vec{v}) = \exp\left[-\frac{1}{2}\left(g(\vec{u},\vec{q})-\vec{y}_{\rm obs}\right)^\prime \vec{R}^{-1}\left(g(\vec{u},\vec{q})-\vec{y}_{\rm obs}\right)\right].
\end{equation}

\subsection{Tempering}
An SMC method suffers when the likelihood $h(\vec{v})$ Eq.~\eqref{eq:weights} is peaked, 
which could be due to very accurate data, amount of data or when the prior poorly approximates the posterior.
A tempering iterative approach tackles this problem by introducing temperatures  $0=\phi_0 < \dots < \phi_T = 1$ 
and corresponding bridging likelihoods $h(\vec{v})^{(\phi_{t}-\phi_{t-1})}$ for $t=1,\dots,T$.
A tempering parameter $\phi_t$ is typically chosen based on effective ensemble size 
\begin{equation}\label{eq:ESS}
 {\rm ESS} := \frac{\left(\sum_{i=1}^{M}w_i\right)^2}{\sum_{i=1}^{M}w_i^2},
\end{equation}
such that ESS does not drop below a certain threshold $1<M_{\rm thresh}<M$.
In order to avoid filter degeneracy, each tempering iteration $t$ needs to be supplied with resampling.\par

\subsection{Deterministic resampling}
Resampling is typically performed by a stochastic approach, which introduces an additional error.  
In TETPF a tempering iteration $t$ is accompanied by a 
{\it deterministic} resampling based on optimal transportation. 
The optimal transport $\vec{S}$ is an $M\times M$ matrix with $s_{ij}$ that satisfy 
\begin{equation}\label{eq:OT1}
s_{ij}\geq 0,\qquad \sum_{i=1}^M s_{ij} = \frac{1}{M}, \qquad   \sum_{j=1}^M s_{ij} = \frac{h\left(\vec{v}^{(t)}_i\right)^{(\phi_{t}-\phi_{t-1})}}{\sum\limits_{j=1}^M h\left(\vec{v}^{(t)}_j\right)^{(\phi_{t}-\phi_{t-1})}},
\end{equation}
and minimizes the cost function 
\begin{equation}\label{eq:OT2}
 \sum_{i,j=1}^{M}s_{ij}\norm{\vec{v}^{(t)}_i-\vec{v}^{(t)}_j}^2
\end{equation}
This gives rise to a resampling with replacement and a stochastic transport matrix $\vec{S}$. 
In order to have a deterministic optimal transformation the following proposal is adopted 
\begin{equation}\label{eq:OT3}
 \tilde{\vec{v}}_j = M\sum_{i=1}^{M}\vec{v}^{(t)}_i \tilde{s}_{ij} \quad  \mbox{for} \quad  j=1,\dots,M,
\end{equation}
where $\tilde{s}_{ij}$ is a solution to the optimization problem Eq.~\eqref{eq:OT1}--\eqref{eq:OT2}. 
To solve the linear transport problem Eq.~\eqref{eq:OT1}--\eqref{eq:OT2}, we use $FastEMD$ algorithm of~\cite{PeWe09}.
Its computational complexity is of order $M^2\ln M$, and the algorithm is available as a $MATLAB$ and a $Python$ subroutine.

\subsubsection{Localization}  
 Ensemble Transform Particle Filter as any particle filter does not have assumption about the posterior.
 Therefore it still demands a large ensemble. For high-dimensional problems this is computationally unfeasible.
 Hence one has to decrease the number of degrees of freedom, i.e. by distance-based localization
 of~\cite{ReCo15,RuDu18} abbreviated here LETPF.
  
 Assume we have a numerical grid of $N\times N$ size with grid cells denoted by $X^l$ for $l=1,\hdots,N^2$. 
 Assume that the uncertain parameter $\vec{u}$ is not grid-based.
 We assume, however, that there exists a matrix $\mathcal{A}$
 such that $\log(\vec{k})=\mathcal{A}\vec{u}$ is grid-based, thus $\log(k^l) = \log[k(X^l)]$.
Then for the local update of an uncertain parameter $\log(k^l)$ we  introduce a diagonal matrix 
$\hat{\textbf{C}}^{l} \in \Reals^{\kappa \times \kappa}$ in the observation space with an element 
 \begin{equation}\label{eq:rho_matrix}
 (\hat{C}^{l})_{\ell,\ell} =\rho \left(\frac{||{{X}^l -{r}^\ell}||}{r^\textrm{loc}}\right)\quad\mbox{for}\quad \ell=1,\dots\kappa.
 \end{equation}
 Here ${r}^\ell$ denotes the location of the observation, $r^\textrm{loc}$ is a localization radius and $\rho(\cdot)$ is a taper function, 
 such as Gaspari-Cohn function by~\cite{GaCo99}
 \begin{equation}\label{eq:GC}
 \rho(r) = 
 \begin{cases}
 1-\frac{5}{3}r^2 + \frac{5}{8} r^3+\frac{1}{2}r^4-\frac{1}{4}r^5, &  \quad 0 \leq r \leq 1, \\
 -\frac{2}{3}r^{-1}+4-5r+\frac{5}{3}r^2+\frac{5}{8}r^3-\frac{1}{2}r^4+\frac{1}{12}r^5, &  \quad 1 \leq r \leq 2, \\
 0, & \quad 2 \leq r.
 \end{cases}  
  \end{equation}
LETPF modifies the likelihood Eq.~\eqref{eq:weights} as following
 \begin{equation}\label{eq:likeloc}
 h^l(\vec{v}) = \exp\left[-\frac{1}{2}\left(g(\vec{u},\vec{q}) - \vec{y}_{\text{obs}}\right)^\prime ( \hat{\textbf{C}}^{l}\textbf{R}^{-1}) \left(g(\vec{u},\vec{q}) - \vec{y}_{\text{obs}}\right)\right],
  \end{equation}
 where $\hat{\textbf{C}}^{l}$ is the diagonal matrix given by Eq.~\eqref{eq:rho_matrix}. 
 Then the optimal transport $\vec{S}^{l}$ is an $M\times M$ matrix with entries $s^l_{ij}$ that satisfy 
\begin{equation}\label{eq:OT1Da1}
 s^l_{ij}\geq 0,\qquad \sum_{i=1}^M s^l_{ij} = \frac{1}{M}, \qquad   \sum_{j=1}^M s^l_{ij} =  \frac{h^l\left(\vec{v}_i^{(t)}\right)^{(\phi_{t}-\phi_{t-1})}}{\sum\limits_{j=1}^M h^l\left(\vec{v}_j^{(t)}\right)^{(\phi_{t}-\phi_{t-1})}},
\end{equation}
and minimizes the cost function 
 \begin{equation}\label{eq:OT1Da2}
 \sum_{i,j=1}^{M} s^l_{ij} \left[\log\left(k_{i}^{l,(t)}\right)-\log\left(k_{j}^{l,(t)}\right)\right]^2.
\end{equation}
 The estimated parameter $\log(\tilde{k}^l)$ is given by
 \begin{equation}\label{eq:OT1Da3}
 \log(\tilde{k}_j^l) = M\sum_{i=1}^{M} \tilde{s}^l_{ij}\log\left(k_i^{l,(t)}\right)\quad \mbox{for}\quad j = 1,\hdots,M,
 \end{equation}
 where $\tilde{s}^l_{ij}$ is is a solution to the optimization problem Eq.~\eqref{eq:OT1Da1}--\eqref{eq:OT1Da2}.
 Then the estimated model parameter is $\tilde{\vec{u}} =  \mathcal{A}^{-1}\log(\tilde{\vec{k}})$.
 
 We note that localization reduces LETPF to a univariate transport problem.
The univariate linear transport problem is solved by sorting the ensemble members~\cite{ReCo15}.
Update of the uncertain grid-based parameters $\log(k^l)$ could be performed in parallel for each grid cell $l=1,\hdots,N^2$. 
Computational complexity of the sorting algorithm per grid cell is $M\ln M$.

To estimate a scalar model error $q$, we solve the optimal transport problem
  with $s^{\rm G}_{ij}$ that satisfy 
\begin{equation}\label{eq:OT1Db1}
 s^{\rm G}_{ij}\geq 0, \qquad \sum_{i=1}^M s^{\rm G}_{ij} = \frac{1}{M}, \qquad   \sum_{j=1}^M s^{\rm G}_{ij} =  \frac{h\left(\left[\tilde{\vec{u}}_i\ q_i^{(t)}\right]^\prime\right)^{(\phi_{t}-\phi_{t-1})}}{\sum\limits_{j=1}^M h\left(\left[\tilde{\vec{u}}_j\ q_j^{(t)}\right]^\prime\right)^{(\phi_{t}-\phi_{t-1})}},
\end{equation}
and minimize the cost function 
 \begin{equation}\label{eq:OT1Db2}
 \sum_{i,j=1}^{M} s^{\rm G}_{ij}\left(\norm{\tilde{\vec{u}}_i - \tilde{\vec{u}}_j}^2 + \left(q_{i}^{(t)}-q_{j}^{(t)}\right)^2\right).
\end{equation}
 The estimated parameter $\tilde{q}$ is given by
 \begin{equation}\label{eq:OT1Db3}
 \tilde{q}_j = M\sum_{i=1}^{M} \tilde{s}^{\rm G}_{ij}q_i^{(t)} ,   \quad j = 1,\hdots,M,
 \end{equation}
 where $\tilde{s}^{\rm G}_{ij}$ is is a solution to the optimization problem Eq.~\eqref{eq:OT1Db1}--~\eqref{eq:OT1Db2}. 
 Finally, we set $\tilde{\vec{v}} = [\tilde{\vec{u}}\ \tilde{q}]^\prime$.
 
 
\subsection{Mutation}
The advantage of deterministic resampling is that it reduces sampling noise. 
The disadvantage of deterministic resampling is that for a deterministic and non-chaotic system the filter collapse is 
unavoidable unless particle mutation is introduced. The mutation is performed over an index $1<\tau<\tau_{\rm max}$
with prescribed $\tau_{\rm max}$. 
At the first inner iteration $\tau = 1$ we assign $\vec{v} = \tilde{\vec{v}}$.

We denote by $v_i^\ell$ a component of $\vec{v}_i$, where $1 \leq \ell \leq n+m$. 
If $v_i^\ell$ has a Gaussian prior, then we use the preconditioned Crank-Nicolson pcn-MCMC method from~\cite{CRSW13} 
\begin{equation}
  v_i^{\ell, \rm prop} = \sqrt{1-\beta^2} v_i^\ell + \beta\xi_i \quad \mbox{for} \quad  i=1,\dots,M, \label{eq:prop}\\
\end{equation}
where $\{\xi_i\}_{i=1}^M$ is from normal distribution. 
If $v_i^\ell$ has a uniform prior $U[a,\ b]$, then we use
\begin{equation}
  v_i^{\ell, \rm prop} =  v_i^\ell + \xi_i \quad \mbox{for} \quad  i=1,\dots,M, \label{eq:propU}\\
\end{equation}
where $\xi_i \sim U[a-b,\ b-a]$, and we project $v_i^{\ell, \rm prop}$ to $[a,\ b]$ when necessary.
The proposal Eq.~\eqref{eq:prop}-- \eqref{eq:propU}
is accepted 
\begin{equation}
  \vec{v} = \vec{v}^{\rm prop} \quad \mbox{with the probability} \quad
 \min\left\{1,\frac{h(\vec{v}^{\rm prop})^{\phi_{t+1}}}{h(\vec{v})^{\phi_{t+1}}} \right\}, \label{eq:acc}
\end{equation}
and the inner iteration $\tau$ is increased by one. 
The mutation Eq.~\eqref{eq:prop}--\eqref{eq:acc} is repeated until $\tau = \tau_{\rm max}$, then
we assign $\vec{v}^{(t+1)} = \vec{v}$. 

After that, next tempering iteration proceeds by computing the weights Eq.~\eqref{eq:weights},
new temperature $\phi$ based on Eq.~\eqref{eq:ESS} ${\rm ESS}\geq M_{\rm thresh}$, performing deterministic resampling either 
by Eq.~\eqref{eq:OT1}--\eqref{eq:OT3} for the non-localized method or by Eq.~\eqref{eq:OT1Da1}--\eqref{eq:OT1Db3} for the localized method,
and concluding by mutation Eq.~\eqref{eq:prop}--\eqref{eq:acc} for $\tau_{\rm max}$ iterations. 
The algorithms stops when the final temperature $\phi$
reaches one. It should be noted that the final tempering iteration $T$ is not predefined but found on the fly.
TETPF demands $T M(\tau_{\rm max}+1)$ evaluations of the model $g$,
and TLETPF demands $T M(\tau_{\rm max}+2)$ evaluations of the model $g$.


\section{Regularized Ensemble Kalman Filter}\label{sec:REnKF}
REnKF is based on the Ensemble Kalman Filter with perturbed observations
\[
  \vec{y}_i^\eta =  \vec{y}_{\rm obs} + \eta_i \quad \mbox{for}\quad  i=1,\dots,M,
\]
where $\eta_i \sim \mathcal{N}(\vec{0},\vec{R})$ with $\vec{R}$ being a known covariance matrix of the observation noise.
We define an $M$-dimensional vector with all elements equal to 1 as $\vec{1}_M$. 
REnKF solves the following set of equations for $t=0,\dots,T-1$ with $\vec{v}^{(0)}$ being an initial ensemble
\begin{align}
    \vec{B}^{\rm gg} &= \frac{1}{M-1}\left(g(\vec{u}^{(t)},\vec{q}^{(t)})-\overbar{g(\vec{u}^{(t)},\vec{q}^{(t)})}\vec{1}_M^\prime\right)\left(g(\vec{u}^{(t)},\vec{q}^{(t)})-\overbar{g(\vec{u}^{(t)},\vec{q}^{(t)})}\vec{1}_M^\prime\right)^\prime,\nonumber\\
    \vec{B}^{\rm vg} &= \frac{1}{M-1}\left(\vec{v}^{(t)}-\overbar{\vec{v}^{(t)}}\vec{1}_M^\prime\right)\left(g(\vec{u}^{(t)},\vec{q}^{(t)})-\overbar{g(\vec{u}^{(t)},\vec{q}^{(t)})}\vec{1}_M^\prime\right)^\prime,\nonumber\\
    \vec{v}^{(t+1)}_i &= \vec{v}^{(t)}_i + \vec{B}^{\rm vg}\left(\vec{B}^{\rm gg} + \mu^{(t)} \vec{R}\right)^{-1}\left(\vec{y}_i^\eta - g(\vec{u}^{(t)}_i,\vec{q}^{(t)}_i)\right) \quad\mbox{for}\quad i=1,\dots,M.\label{eq:EnKF}
\end{align}
The regularized parameter $\mu^{(t)}$ is chosen such that
\begin{equation}\label{eq:dp}
 \mu^{(t)} \norm{\ \vec{R}^{1/2} \left(\vec{B}^{\rm gg} + \mu^{(t)} \vec{R}\right)^{-1}\left(\vec{y}_{\rm obs} - \overbar{g(\vec{u}^{(t)},\vec{q}^{(t)})}\right)}
 \geq \Omega \norm{ \vec{R}^{-1/2} \left(\vec{y}_{\rm obs} - \overbar{g(\vec{u}^{(t)},\vec{q}^{(t)})}\right)}
\end{equation}
for predefined $\Omega \in(0, 1)$. 
This is achieved by the bysection method $\mu^{\tau+1} = 2^\tau\mu^0$ for $\tau=0,\dots,\tau_{\rm max}$
and an initial guess $\mu^0$.
We assign $\mu^{(t)} = \mu^{\tau_{\rm max}}$,  
where $\tau_{\rm max}$ is the first integer for which Eq.~\eqref{eq:dp} holds.

Finally, REnKF is stopped based on discrepancy principle, namely when
\begin{equation*}\label{eq:stop}
 \norm{ \vec{R}^{-1/2} \left(\vec{y}_{\rm obs} - \overbar{g(\vec{u}^{(t)},\vec{q}^{(t)})}\right)}
 \leq 1/\Omega \norm{\vec{R}^{-1/2} \eta}
\end{equation*}
with $\eta$ being the observation noise. 
The rule of thumb is to choose $\Omega\in(0.5, 1)$, and we choose $\Omega = 0.7$ for all the numerical experiments.
REnKF demands $T M +1$ evaluations of the model $g$.

\subsection{Localization}  
Covariance-based localization~\cite{HaWh01,HoMi01} can be applied to an Ensemble Kalman filter in order to
remove spurious correlations due to a small ensemble size. 
We assume again having a numerical grid of $N\times N$ size with grid cells denoted by $X^l$ for $l=1,\hdots,N^2$. 
 Assume that the uncertain parameter $\vec{u}$ is not grid-based.
 We assume, however, that there exists a matrix $\mathcal{A}$
 such that $\log(\vec{k})=\mathcal{A}\vec{u}$ is grid-based, thus $\log(k^l) = \log[k(X^l)]$.
Then Eq.~\eqref{eq:EnKF} for a localized EnKF, denoted here LEnKF is rewritten as
\begin{align*}
  \log\left({\vec{k}}^{(t+1)}_i\right) &= \log\left(\vec{k}^{(t)}_i\right) + \hat{\mathcal{C}}\circ\vec{B}^{\rm \log(k)g}\left(\vec{B}^{\rm gg} + \mu^{(t)} \vec{R}\right)^{-1}\left(\vec{y}_i^\eta - g(\vec{u}^{(t)}_i,\vec{q}^{(t)}_i)\right) \quad\mbox{for}\quad i=1,\dots,M,\\
  \eta^{(t+1)}_i &= \eta^{(t)}_i + \vec{B}^{\rm \eta g}\left(\vec{B}^{\rm gg} + \mu^{(t)} \vec{R}\right)^{-1}\left(\vec{y}_i^\eta - g(\vec{u}^{(t)}_i,\vec{q}^{(t)}_i)\right) \quad\mbox{for}\quad i=1,\dots,M. 
\end{align*}
Here $\circ$ denotes the element-wise product and $\hat{\mathcal{C}}$ is a distance-based correlation matrix, an element of which is 
 \begin{equation}\label{eq:rho_matrixKF}
 \hat{\mathcal{C}}_{l,\ell} =\rho \left(\frac{||{{X}^l -{r}^\ell}||}{r^\textrm{loc}}\right)\quad\mbox{for}\quad l=1,\dots,N^2\quad \mbox{and}\quad  \ell=1,\dots,\kappa,
 \end{equation}
where ${r}^\ell$ denotes the location of the observation, $r^\textrm{loc}$ is a localization radius and $\rho$ is given by Eq.~\eqref{eq:GC}.
The covariance matrices $\vec{B}^{\rm \log(k)g}$ and $\vec{B}^{\rm \eta g}$ are 
\begin{align*}
    \vec{B}^{\rm \log(k)g} &= \frac{1}{M-1}\left( \log(\vec{k}^{(t)})-\overbar{\log(\vec{k}^{(t)})}\vec{1}_M^\prime\right)\left(g(\vec{u}^{(t)},\vec{q}^{(t)})-\overbar{g(\vec{u}^{(t)},\vec{q}^{(t)})}\vec{1}_M^\prime\right)^\prime,\\
    \vec{B}^{\rm \eta g} &= \frac{1}{M-1}\left(\vec{\eta}^{(t)}-\overbar{\vec{\eta}^{(t)}}\vec{1}_M^\prime\right)\left(g(\vec{u}^{(t)},\vec{q}^{(t)})-\overbar{g(\vec{u}^{(t)},\vec{q}^{(t)})}\vec{1}_M^\prime\right)^\prime.
\end{align*}
RLEnKF also demands $T M +1$ evaluations of the model $g$, as REnKF.

\section{Numerical experiments}\label{sec:NE}
In this Section we apply T(L)ETPF to nonlinear problems of
different dimensionality and compare it to R(L)EnKF.
A first problem is low-dimensional and a second one is high-dimensional.
In the first problem, model error is multiplicative, and the model consists of two scalar-valued functions.
This problem is not grid-based, therefore localization cannot be applied. 
In the second problem, the source of model error is uncertain boundary conditions, and the model itself is a steady-state single-phase Darcy flow model.
This problem is grid-based, therefore localization can be applied. 

\subsection{Multiplicative model error}
We consider a test case of estimating two uncertain parameters with non-Gaussian marginalized posterior.
Consider a model with multiplicative model noise
\[
\vec{g}(\vec{u},\vec{q}) = \begin{bmatrix}
 g(u^1,q^1)\\
 g(u^2,q^2)
\end{bmatrix},  \quad \mbox{where} \quad g(u,q) = q\exp\left[1-\frac{9}{2} \left(u-\frac{2 \pi}{3} \right)^2\right].
\]
The true model is $\vec{G}(\vec{u}) = \vec{q}(\vec{u},\vec{1})$.
Observations are $y^\ell_{\rm obs} = 1.8$ for $\ell=1,2$ and the observation noise is $\mathcal{N}(0,0.001)$.
The prior for uncertain parameters $u^\ell$ is $\mathcal{N}(2.4,1)$ for $\ell=1,2$ 
and for model error $q^\ell$ is $\mathcal{N}(1,0.01)$ for $\ell=1,2$.
For this toy problem the true posterior can be computed directly by 
the Bayes approach Eq.~\eqref{eq:exp}--\eqref{eq:joint} with a large sample of $\vec{u}$ and $\vec{q}$ of size $10^4$ each.
Thus the joint probability density function is computed over a space of dimension $10^4\times 10^4$.

For TETPF, we choose the threshold for ESS to be $M_{\rm thresh} = M/2$.
Parameters $\tau_{\rm max} = 20$ and $\beta = 0.02$ give good mixing and thus
were used for the mutation step.
For both TETPF and REnKF we use ensemble size $M=1000$.
We perform ten numerical experiments to check initial sample sensitivity.
On average TETPF took eight tempering iterations and REnKF ten regularizing iterations.

\begin{figure}[ht]
	\centering
{\includegraphics[scale=0.4]{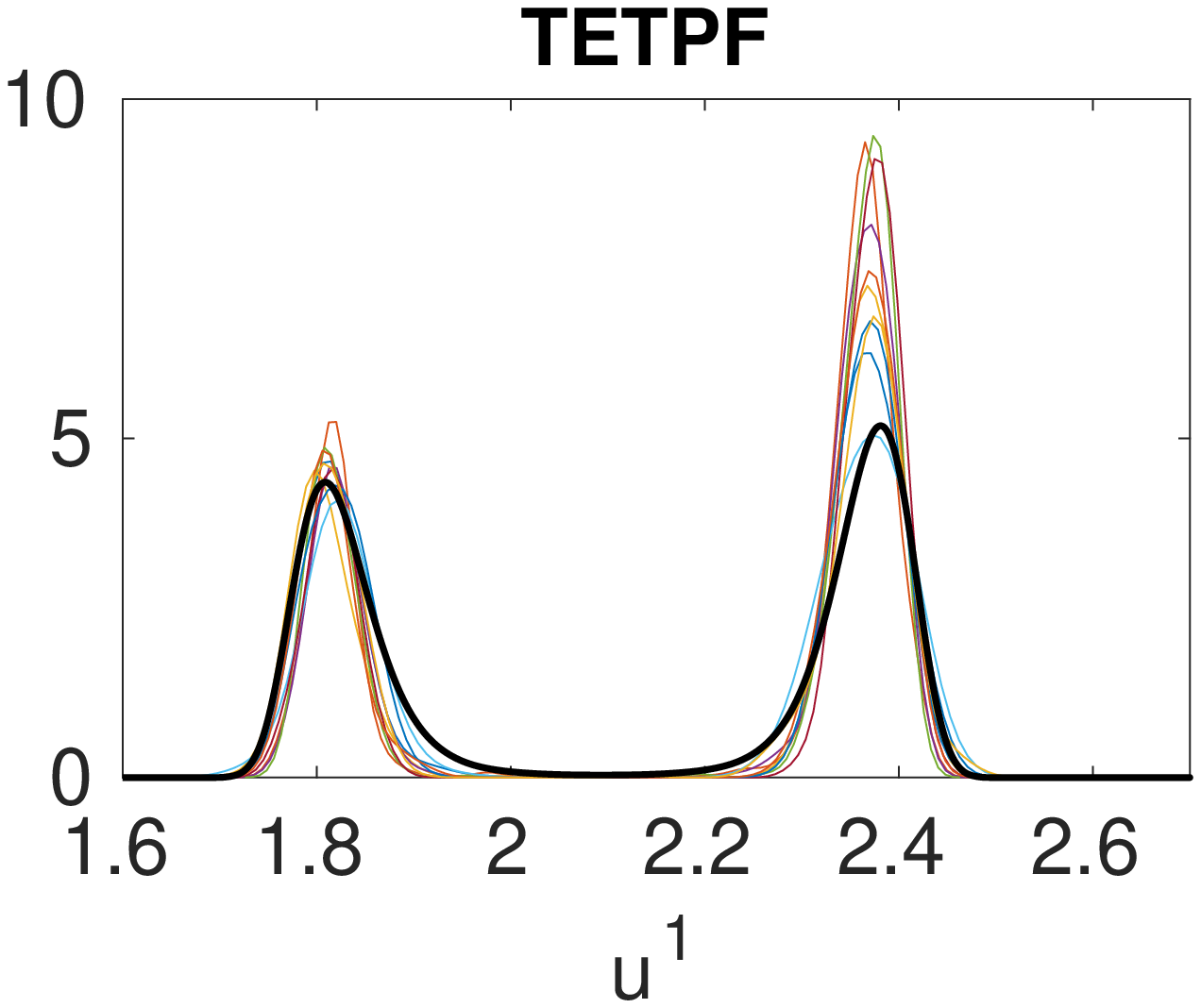}}
{\includegraphics[scale=0.4]{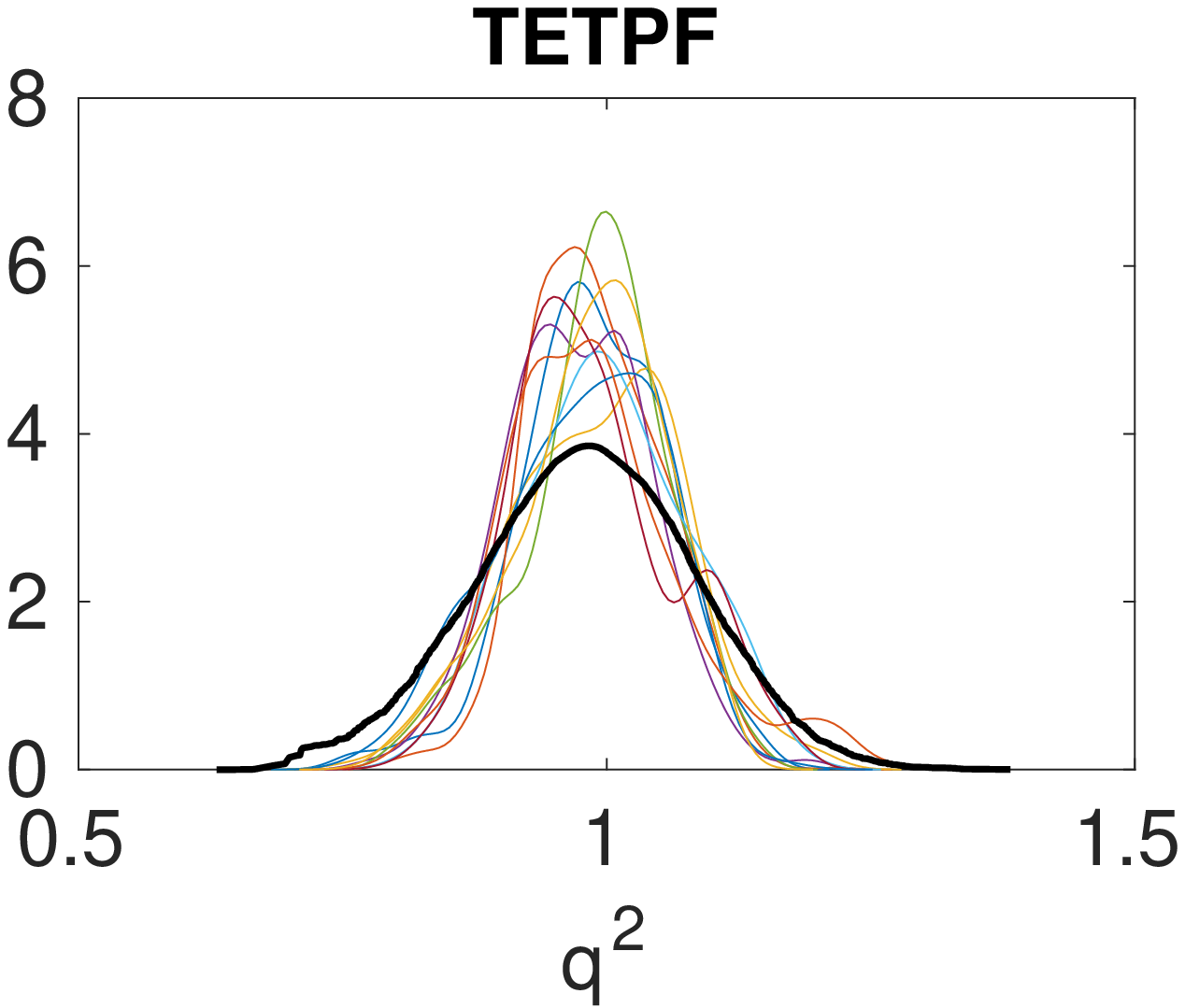}}
{\includegraphics[scale=0.4]{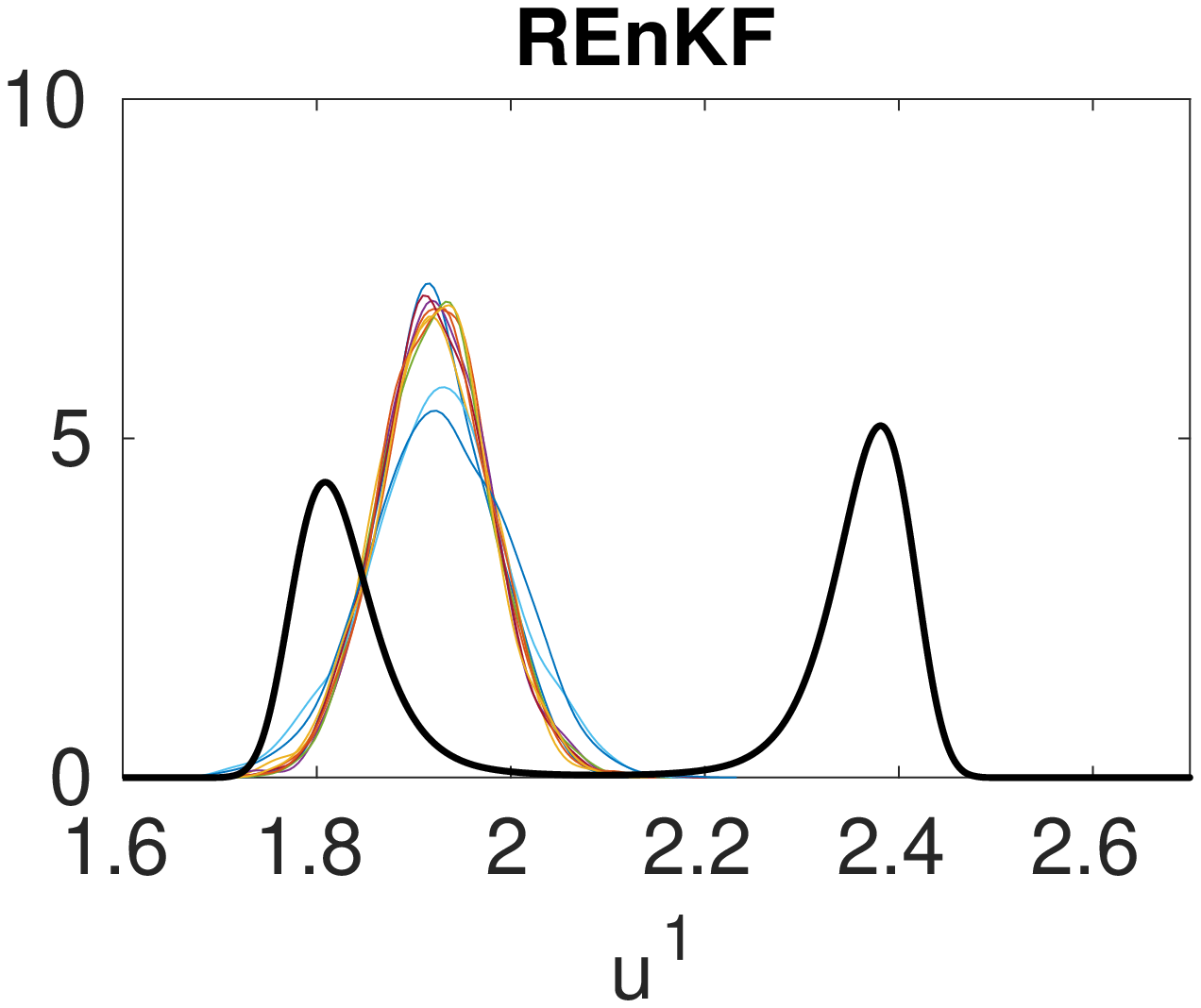}}
{\includegraphics[scale=0.4]{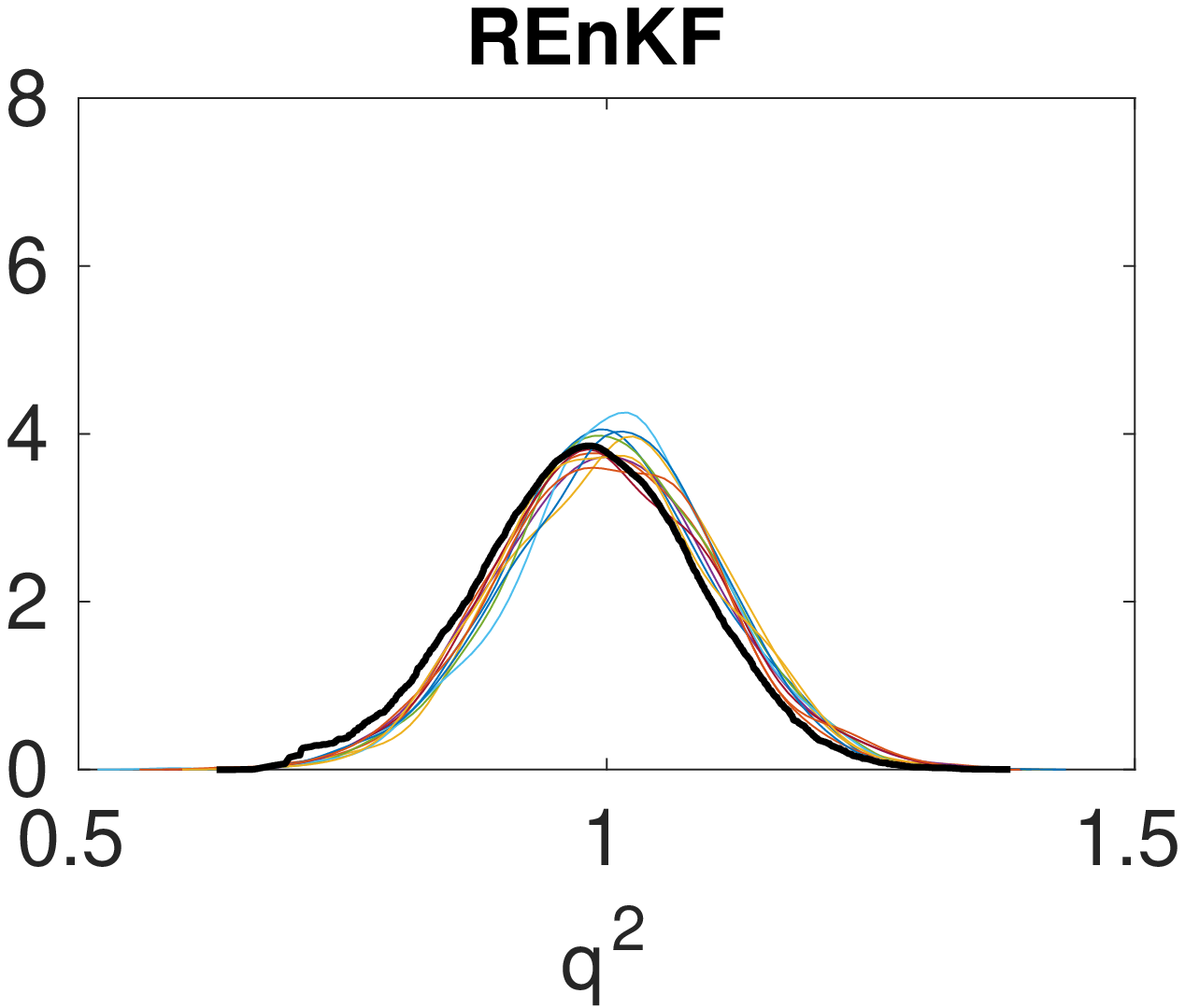}}
	\caption{Posterior of $u^1$ on the left and  of $q^2$ on the right for the toy problem. 
	In black is the true probability, and in color is an SMC approximation with different colors corresponding
	to different simulations. At the top is TETPF and at the bottom is REnKF.}\label{fig:Fig1}
\end{figure}

In Fig.~\ref{fig:Fig1} we plot posteriors for $u^1$ on the left and for $q^2$ on the right ($u^2$ and $q^1$ give similar results
and hence are omitted).
We observe that TETPF gives good approximation of uncertain parameters and these approximations are better than the ones given by REnKF. 
For model error, REnKF is on the contrary robust compared to TETPF. 
Therefore, REnKF fails to estimate non-Gaussian probability for the uncertain parameters, while gives a good approximation for the Gaussian distributed model error.
This is to be expected since an Ensemble Kalman Filter has an assumption of Gaussian probabilities and has been proven to be highly efficient and 
robust for estimating those  probabilities.
TETPF does not have such an assumption. Therefore it has larger error than REnKF due to fully unconstrained optimization problem. 
This leads to a multimodal estimation on the one hand but larger noise than REnKF for a Gaussian probability on the other hand. 
\subsection{Uncertain boundary conditions}\label{sec:Darcy}
We consider a test case of estimating uncertain Gaussian permeability. 
We consider a steady-state single-phase Darcy flow model defined over an aquifer of two-dimensional physical domain  $D = [0,6] \times [0,6]$, which is given by
\begin{equation*}\label{eq:Darcy}
-\nabla \cdot \left[k(x,y)\nabla P(x,y)\right]  =\mathcal{F}(x,y) \quad \mbox{for}\quad  (x,y)\in D,  \\
\end{equation*} 
where $\nabla = (\partial/\partial{x} \ \partial/\partial{y})^\prime$, $\cdot$ is the dot product, $P(x,y)$ is the pressure, 
$k(x,y)$ is the permeability, and the source term is
\[
\mathcal{F}(x,y) =
\begin{cases}
0 \quad &\mbox{for} \quad 0 \leq y \leq 4,\\
137 \quad &\mbox{for} \quad 4<y \leq 5,\\
274 \quad &\mbox{for} \quad 5<y \leq 6.
\end{cases}
\] 
The boundary conditions are a combination of Dirichlet and Neumann boundary conditions 
\begin{equation*}\label{eq:DNbc}
P(x,0) = 100, \quad \frac{\partial P}{\partial x} (6,y) = 0, \quad -k(0,y) \frac{\partial P}{\partial x} (0,y) = 500(1+q), \quad \frac{\partial P}{\partial y} (x,6) = 0,
\end{equation*}
where $q$ denotes model error. The true model $G$ does not have the error in the boundary conditions, while the incorrect model $g$ does. 
Thus we have $G(u) = g(u,0)$,
where $u$ is an uncertain parameter related to permeability.

We implement a cell-centered finite difference method to discretize the domain $D$ into $N\times N$ grid cells $X^l$ of size $\Delta x^2$.
We  solve the true forward model $G$ on a fine grid $N=N_{\rm f}=140$ for the true solution.
Then the synthetic observations are obtained by
\begin{equation*} 
\vec{ y}_{\rm obs} = \vec{ L}(\vec{ P}^{\rm true}) + \eta.              
\end{equation*}
An element of 
$\vec{L}(\vec{P}^{\rm true})$ is a linear functional of pressure, namely
\begin{equation*}
L^\ell(\vec{ P}^{\rm true}) = \frac{1}{2 \pi \sigma ^2} \sum_{l=1}^{N_{\rm f}^2} \exp \left(-\frac{||X^l-r^\ell ||^2}{2 \sigma^2}\right) P^{{\rm true},l} \Delta x_{\rm f}^2 \quad \mbox{for} \quad \ell=1, \dots,\kappa,
\end{equation*}
where $r^\ell$ denotes the location of the observation, $\kappa$ the number of observations, $\Delta x_{\rm f} = 6/N_{\rm f}$, and $\sigma = 0.01$.
Observation noise is denoted by $\eta$ and it is drawn from $\mathcal{N}(\vec{0},\vec{R})$. 
Observation error covariance $\vec{R}$ is known, and we choose it such that 
the norm of the noise is 1\% of the norm of the data.

Both the true permeability and an initial ensemble are drawn from the same prior distribution
as the prior includes knowledge about geological properties. 
We assume log permeability is generated by a random draw from a Gaussian distribution $\mathcal{N}(\log(\vec{5}),\vec{C})$.
Here  $\vec{5}$ is an $N^2$ vector with all elements being 5
and $\vec{C}$ is Whittle-Matern correlation, an element of which is given by
\begin{equation*}
	{C}^{l\ell} = \frac{1} {\Gamma (1)} \frac{d^{l\ell}}{\delta} \Upsilon_1 \left(\frac{d^{l\ell}}{\delta}\right) \quad \mbox{for} \quad  l,\ell=1,\dots,N^2. 
\end{equation*}
Here $d^{l\ell}$ is the distance between two spatial locations, $\delta=0.5$ is the correlation length,
$\Gamma$ is the gamma function, and
$\Upsilon_1$ is the modified Bessel function of the second kind of order one.
We denote by $\lambda$ and $\gamma$ eigenvalues and eigenfunctions of $\vec{C}$, respectively,
then following  Karhunen-Loeve expansion log permeability is
\begin{equation*}\label{eq:KLexp}
\log(k^l) = \log(5)+ \sum_{\ell=1}^{N^2} \sqrt{\lambda^\ell} \gamma^{\ell l}  u^\ell \quad \mbox{for}  \quad l=1,\dots,N^2,
\end{equation*}
where $u^\ell$ is i.i.d. from a normal distribution for $\ell=1,\dots,N^2$.

Therefore the initial parameter $u$ is drawn from $\mathcal{N}(0,1)$, while 
the initial model error $q$ is drawn from a uniform distribution $U[0\ 0.5]$.
We then solve the incorrect forward model $g$ on a coarse grid $N=N_{\rm c}=70$.
The uncertain parameter $\vec{u}$ has the dimension $n=4900$,
which makes the dimension of $\vec{v}$ $n+m=4901$.
We perform 20 numerical experiments with both T(L)ETPF and R(L)EnKF to check initial sample sensitivity.
We conduct numerical experiments with ensemble sizes 100 and 1000.
We compare the methods to a pcn-MCMC method. An MCMC experiment
was conducted using 200 chains with the lengths $10^6$, burn-in period $10^5$, and thinning period $10^3$ each. 
For T(L)ETPF we choose $\tau_{\rm max} = 20$ and $\beta = 0.045$ for mutation, since
it gives good mixing with acceptance rate at the final tempering iteration around 0.2.
We set the threshold for ESS to be $M_{\rm thresh} = M/3$. 

We define the root mean square error (RMSE) of a mean field $\overline{\Xi}=1/M\sum_{i=1}^M \Xi_i$ as
\[
\text{RMSE}(\Xi) =\sqrt[]{\left(\overline{\Xi}-\Xi^\textrm{MCMC}\right)^T \left( \overline{\Xi}-\Xi^\textrm{MCMC}\right)}
\]
for either log permeability $\Xi = \log({\vec{k}})$ or pressure $\Xi = \vec{P}$.
To choose a favoring localization radius, we perform a numerical experiment with $r^{\rm loc}$
ranging from one to six with an increment of one. Then we define the 
favoring localization radius as a localization radius that gives the smallest RMSE in terms of mean log permeability
for that numerical experiment.
For TLETPF the favoring localization radius is $r^{\rm loc}=1$ for both ensemble sizes 100 and 1000.
For RLEnKF the favoring localization radius is $r^{\rm loc}=3$ for both ensemble sizes 100 and 1000.
 
\begin{figure}[ht]
	\centering
{\includegraphics[scale=0.4]{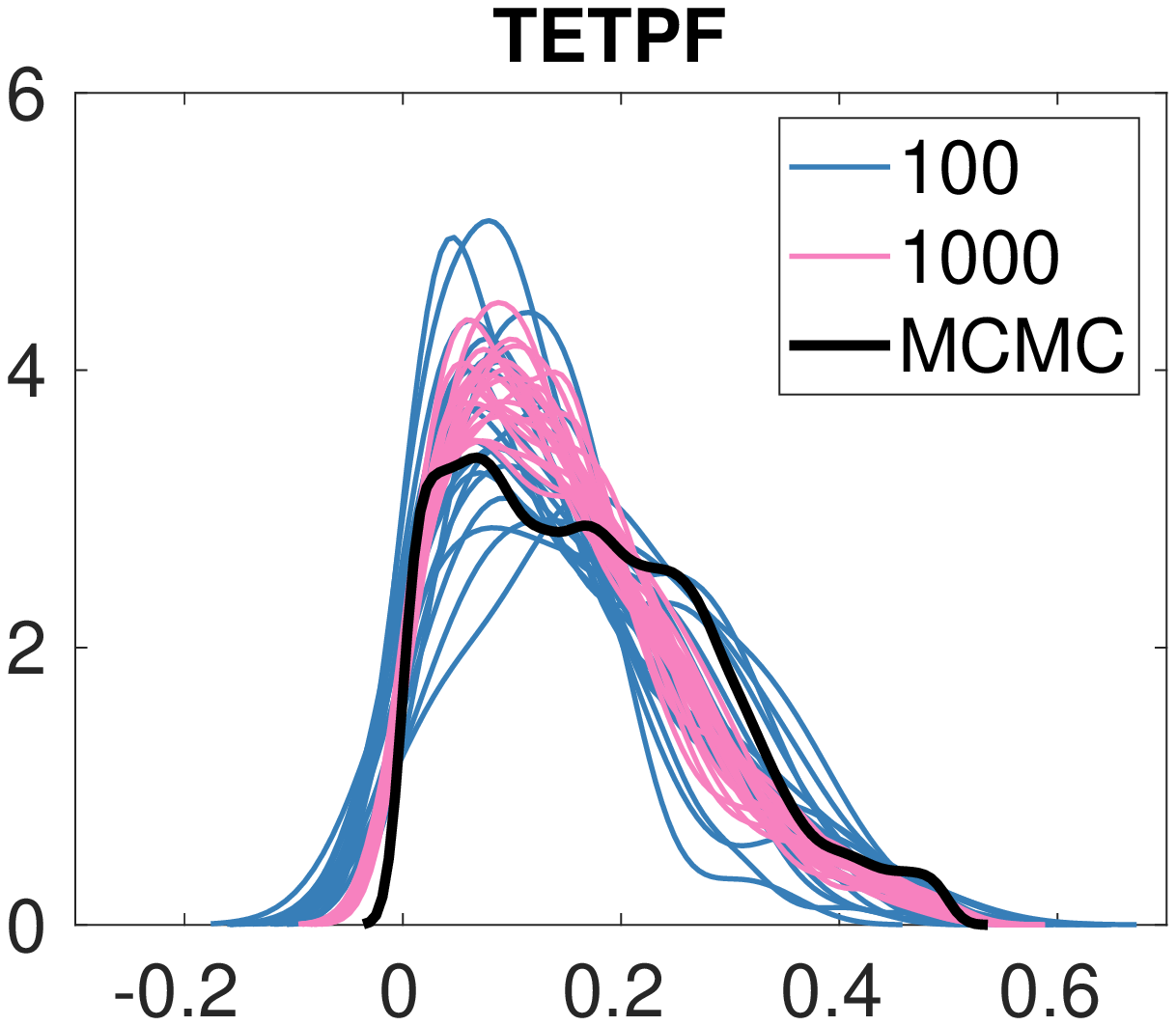}}
{\includegraphics[scale=0.4]{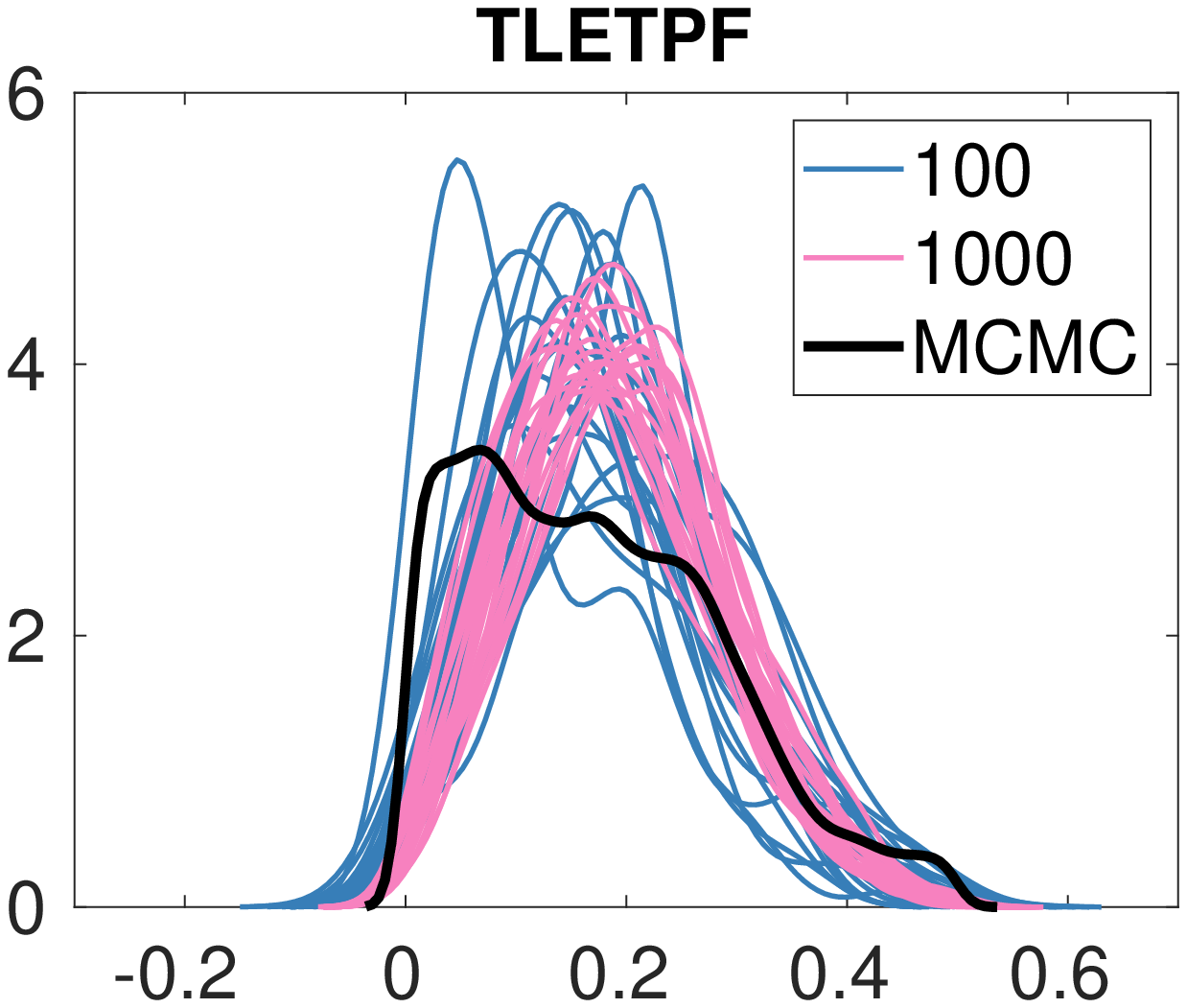}}
{\includegraphics[scale=0.4]{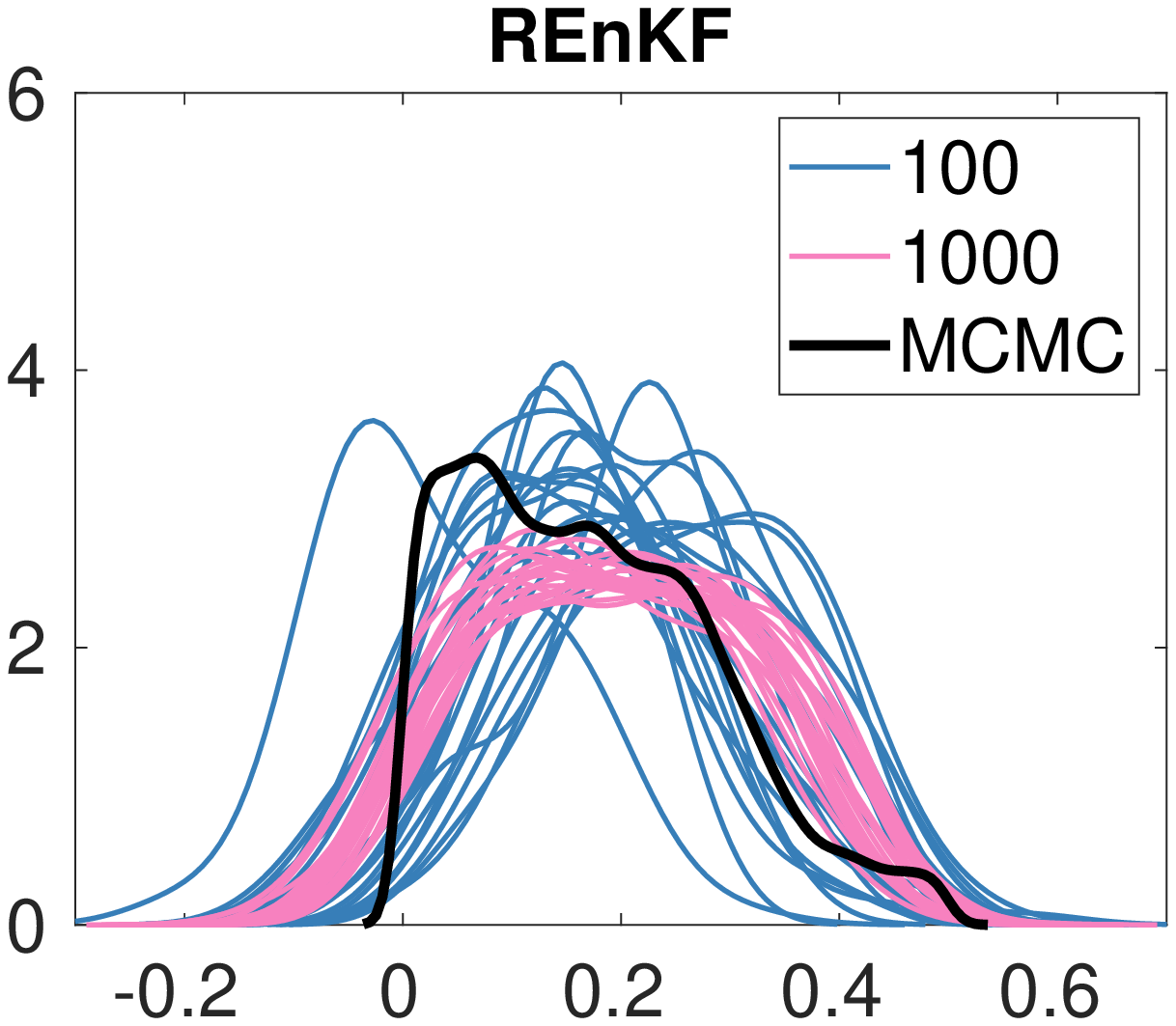}}
{\includegraphics[scale=0.4]{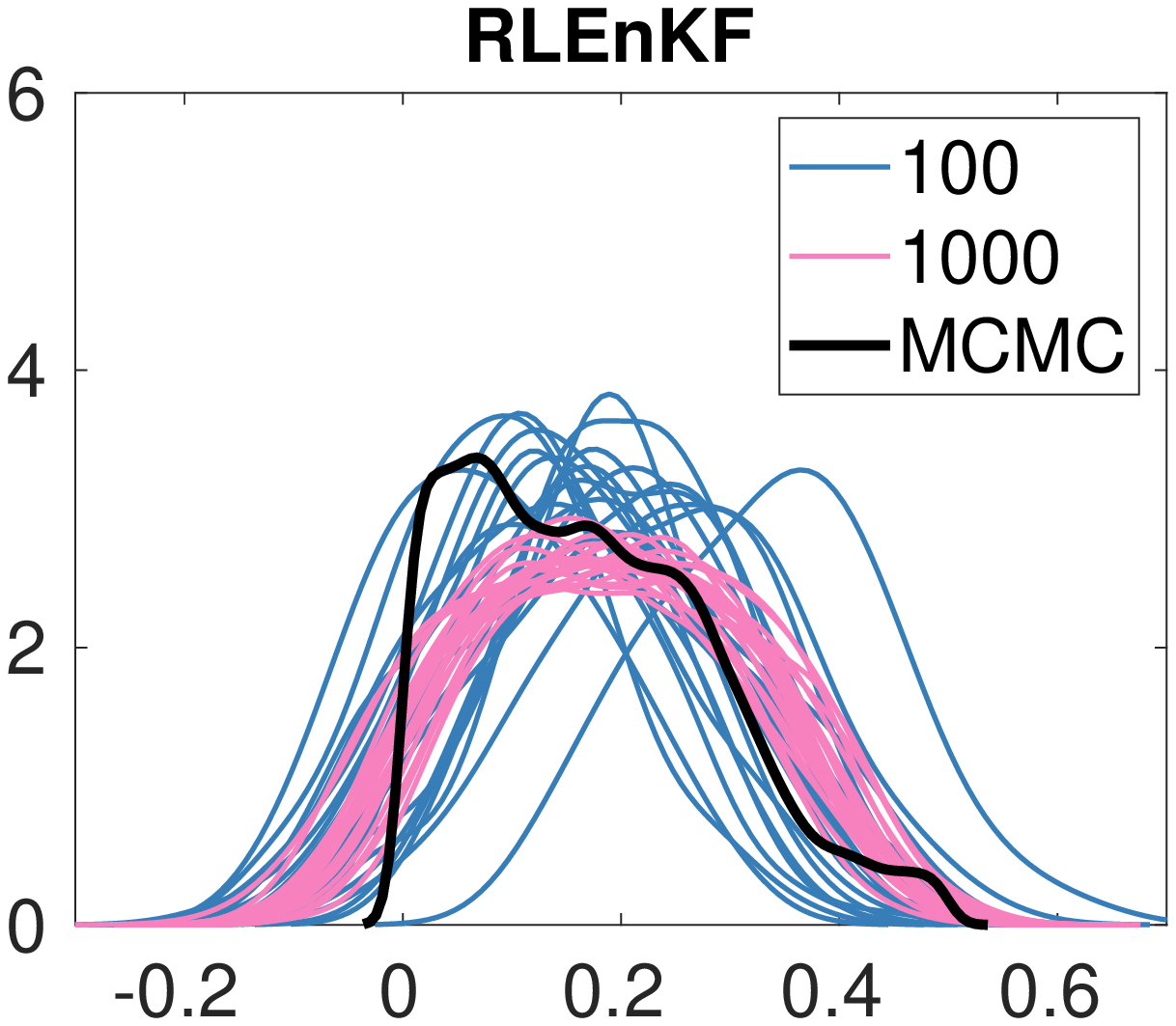}}
\caption{Posterior of model error $q$ for the Darcy flow problem. 
	At the top is T(L)ETPF and at the bottom is R(L)EnKF.
	On the left are the non-localized methods and on the right the localized methods.
	Results for the ensemble size 100 are shown in blue and for 1000 in pink, where one line is for one simulation out of twenty.
	MCMC is shown in black.}\label{fig:Fig2}
\end{figure}

First, we investigate the methods performance with respect to estimation of non-Gaussian model error.  
In Fig.~\ref{fig:Fig2}, we plot the posterior approximations. The MCMC posterior is skewed, while the 
posterior of R(L)EnKF is Gaussian. TETPF gives the best posterior approximation, while TLETPF deteriorates the
results as it was already observed in~\cite{RuDu18}, which makes localization a necessary but still evil.

\begin{figure}[ht]
	\centering
{\includegraphics[scale=0.3]{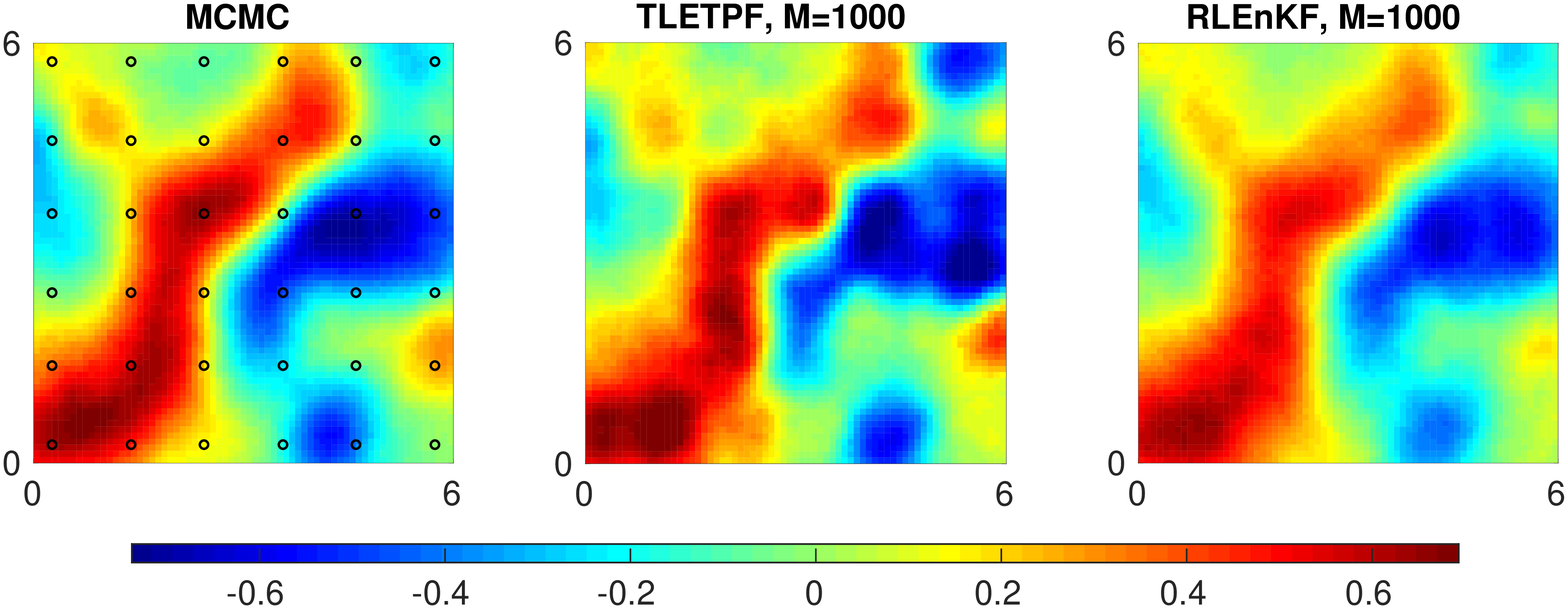}}
{\includegraphics[scale=0.3]{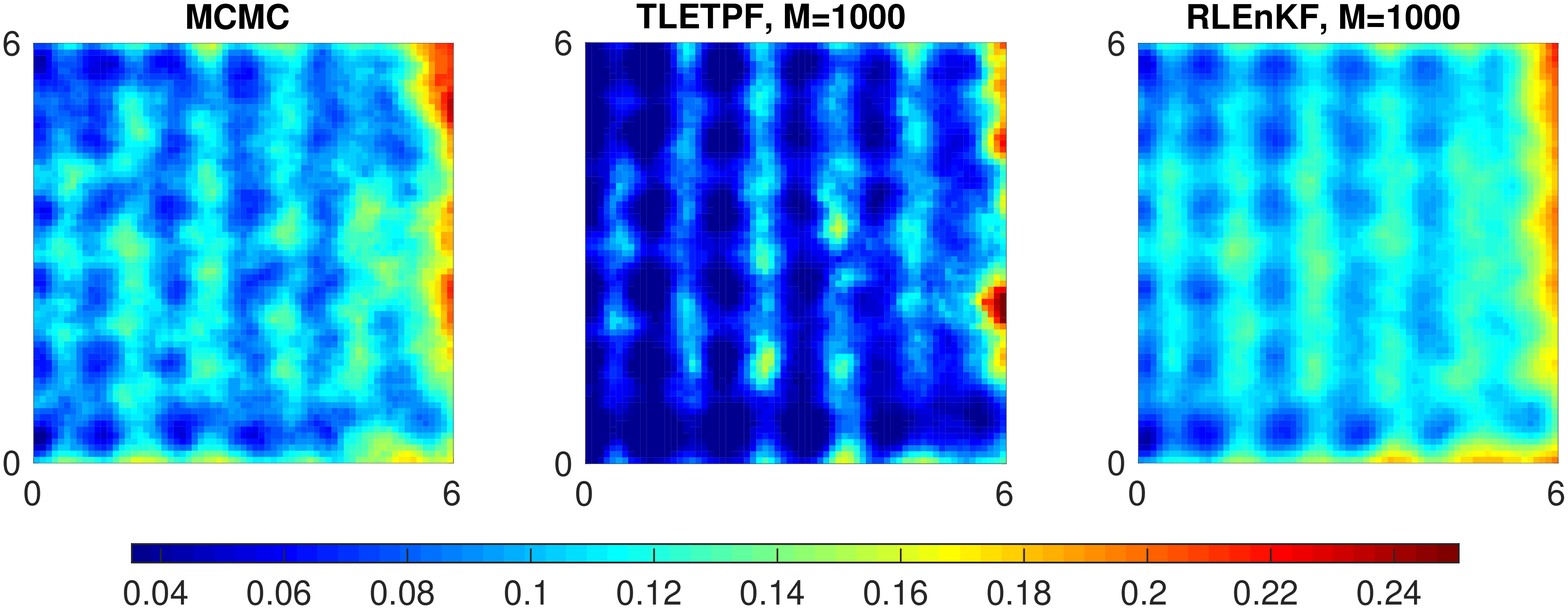}}
\caption{Mean $\log(\vec{k})$ at the top and variance of $\log(\vec{k})$ at the bottom.
MCMC is on the left, with circles for the 
observation locations. A simulation with smallest RMSE of mean log permeability is in the middle for T(L)ETPF and on the right for R(L)EnKF.}\label{fig:Fig3}
\end{figure}

Next, we compare estimations of log permeability. We compute RMSE of mean log permeability. 
A simulation that gives the smallest RMSE is chosen to display results of the mean field and
variance. TLETPF and RLEnKF at ensemble size 1000 give the smallest RMSE.
In Fig~\ref{fig:Fig3}, we plot mean of log permeability at the top
and variance of log permeability at the bottom for MCMC, TLETPF, and RLEnKF.
We observe that both methods give a reasonably good approximation of the MCMC mean log permeability.
The variance, however, is underestimated. 

\begin{figure}[ht]
	\centering
{\includegraphics[scale=0.33]{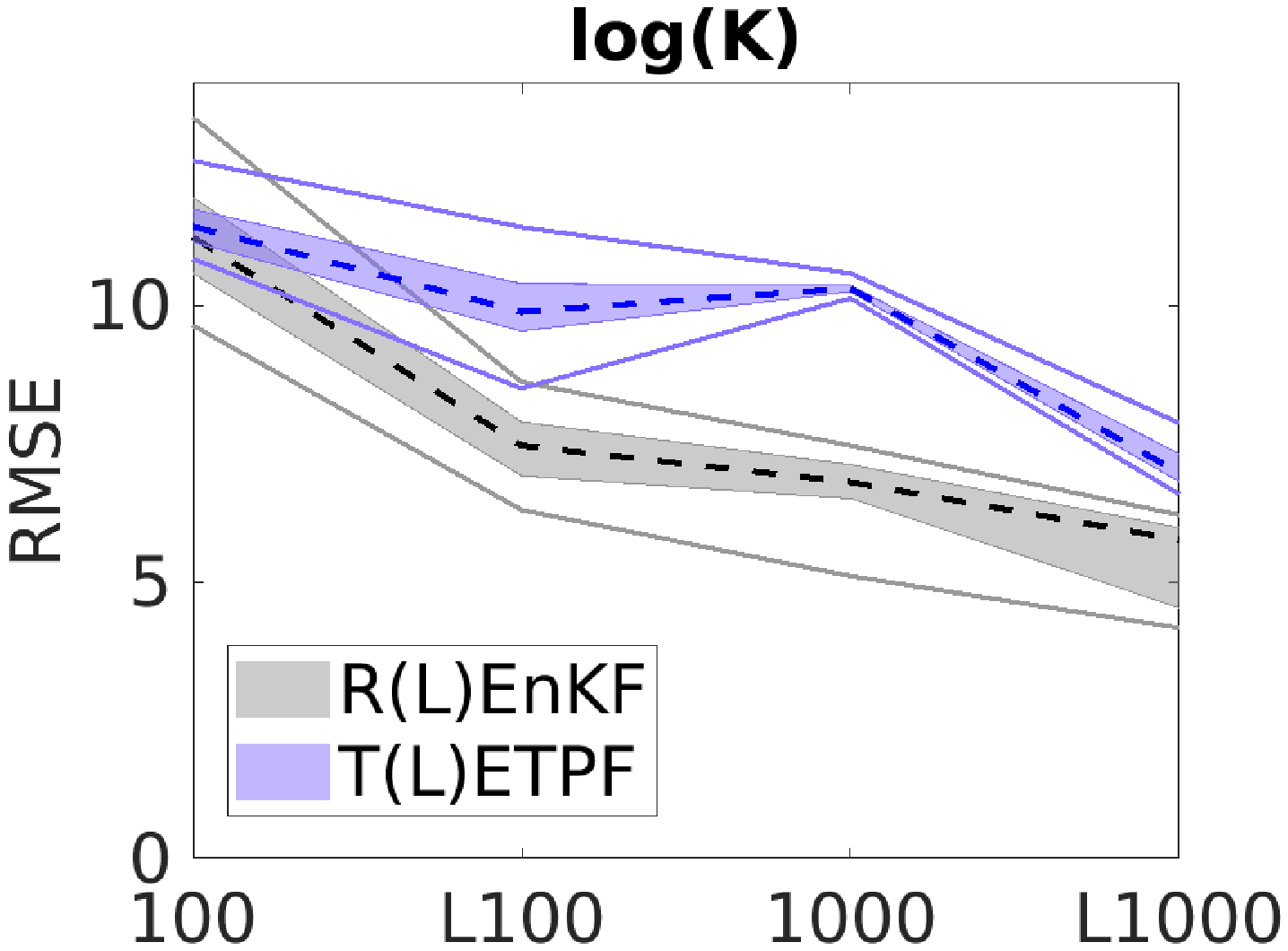}}
{\includegraphics[scale=0.33]{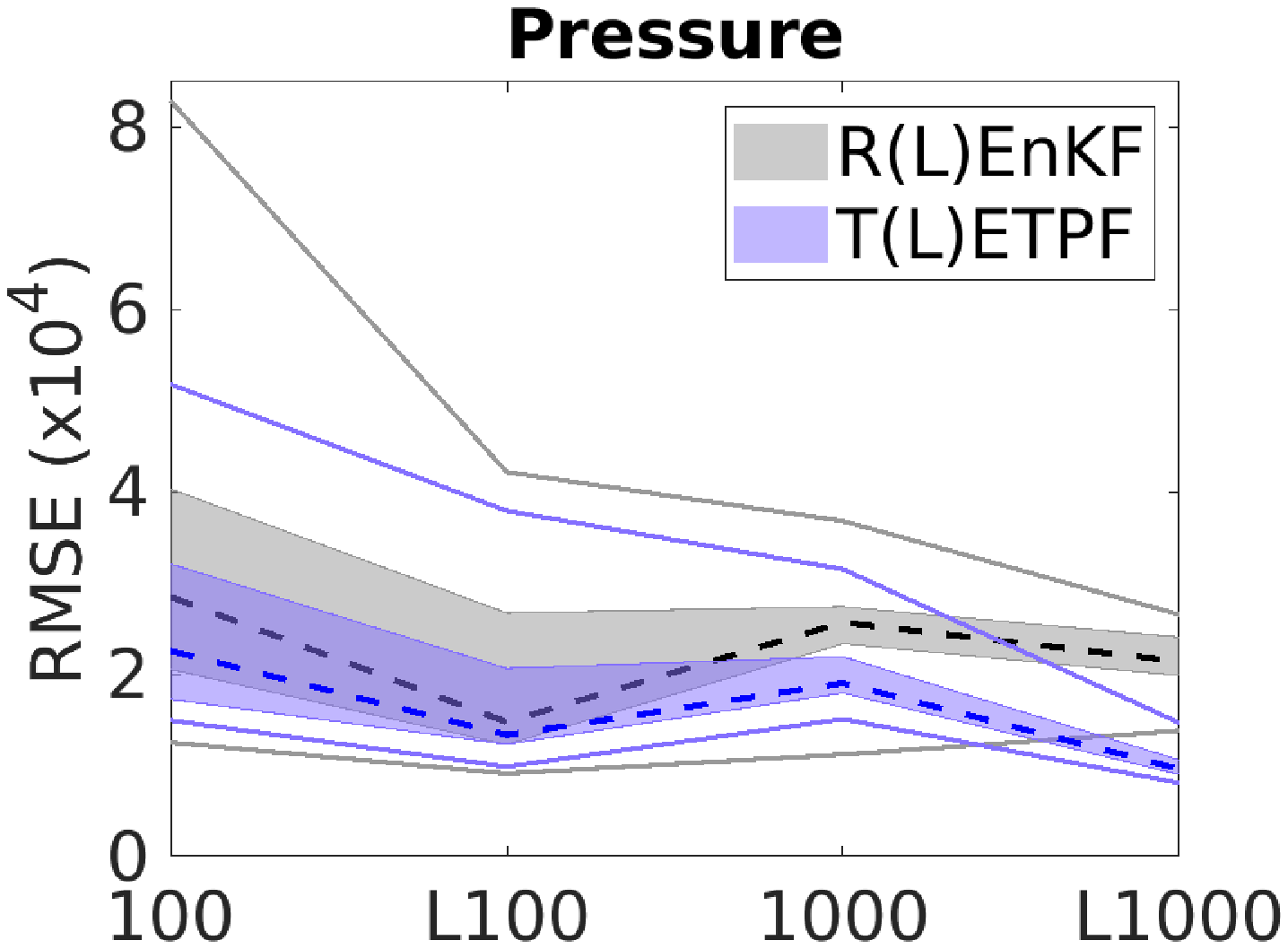}}
\caption{RMSE for $\log(\vec{k})$ is 
on the left. RMSE $\times10^4$ for pressure is on the right. 
A dashed line is for the median,
a shaded area is for 25 and 75 percentile, and a solid line is for 2 and 98 percentile over 20 simulations. 
On the x-axis numbers stand for ensemble sizes and L stands for a localized method.  
R(L)EnKF is shown in gray. T(L)ETPF is shown in blue.}\label{fig:Fig4}
\end{figure}

In Fig.~\ref{fig:Fig4}, we plot RMSE of mean log permeability on the left and of mean pressure on the right at different ensemble sizes
for both T(L)ETPF and R(L)EnKF. 
First, we observe that localization decreases the error, as it was already extensively reported in
the literature. Comparing RMSE of mean log permeability, we see that R(L)EnKF outperforms T(L)ETPF at any ensemble size.
This is again something to be expected, since log permeability is described by a Gaussian process and R(L)EnKF
is excellent at predicting Gaussian probabilities. 
The question that we ask, however, is does this excellent estimation of permeability compensates for a poor estimation of uncertain boundary conditions shown
at the bottom of Fig.~\ref{fig:Fig2}? Comparing RMSE of mean pressure, we observe that it does not, since T(L)ETPF gives smaller RMSE than R(L)EnKF 
at any ensemble size.
Therefore, estimations of pressure {\it are} sensitive to uncertainty in boundary conditions even at one boundary. 
Moreover, even though T(L)ETPF gives worse estimations of permeability than R(L)EnKF, 
by correctly estimating model error T(L)ETPF gives more accurate estimations of pressure.

In Tab.~\ref{tab:iter}, we show the number of iterations 
a method takes on average. First, we would like to remark that REnKF takes less iterations than RLEnKF,
and as the localization radius decreases the number of iterations increases (not shown).
For particle filter, the opposite holds: TETPF takes more iterations than TLETPF. 
This is due to a fundamental difference in localization between RLEnKF and TLETPF.
\begin{table}[h]
\caption{Number of iterations for T(L)ETPF and R(L)EnKF at different ensemble sizes M.\label{tab:iter}}
\begin{tabular}{lllll}
\hline\noalign{\smallskip}
M & TETPF & TLETPF   & REnKF &RLEnKF  \\
\noalign{\smallskip}\hline\noalign{\smallskip}
	 100 & 62 & 35 & 11 & 15\\ 
	 1000 & 65 & 35 & 10 & 13\\	
	 6000 & -- & -- & -- & 11\\ 
	 7700 & -- & -- & 10 & --\\	
\noalign{\smallskip}\hline
\end{tabular}
\end{table}

Let us assume the observation locations $r^\ell$ coincide with some grid cells $X^l$.
We consider a limiting case when the taper function $\rho$ is a Dirac delta function.
Then the localization covariance matrix $\hat{\mathcal{C}}$ from Eq.~\eqref{eq:rho_matrixKF} is a matrix
with all zeros but ones at the observation locations. This means not only that we take one observation 
into account but also that we do not have correlations between parameters due to the element-wise product $ \hat{\mathcal{C}}\circ\vec{B}^{\rm \log(k)g}$.
Therefore, the smaller localization radius is the noisier an RLEnKF approximation becomes.
This makes the optimization problem harder to solve, and thus the convergence of the regularizing parameter $\mu^{(t)}$
is slower, which in turn results in iteration increase.

TLETPF, on the contrary, does not explore any correlations between parameters independent of localization radius, 
since Eq.~\eqref{eq:OT1Da1}--\eqref{eq:OT1Da3} is a univariate optimization problem.
In the limiting case described above the likelihood Eq.~\eqref{eq:likeloc} is less 
picked compared to the likelihood Eq.~\eqref{eq:weights} of TETPF due to a fewer observations taken into account.
This means that the temperature $\phi^{(t)}$ converges faster, which in turn results in iteration decrease.

\begin{figure}[ht]
	\centering
{\includegraphics[scale=0.33]{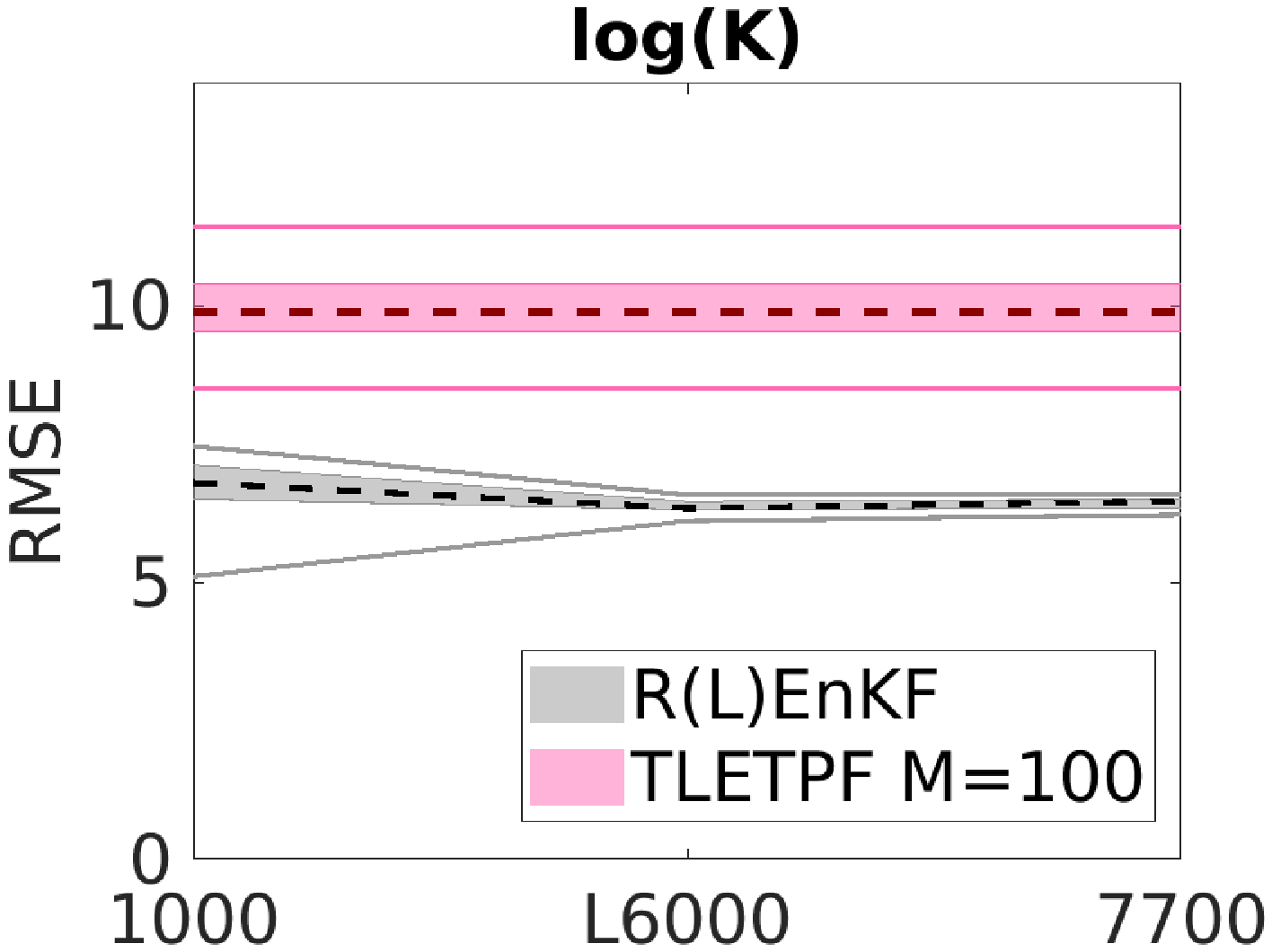}}
{\includegraphics[scale=0.33]{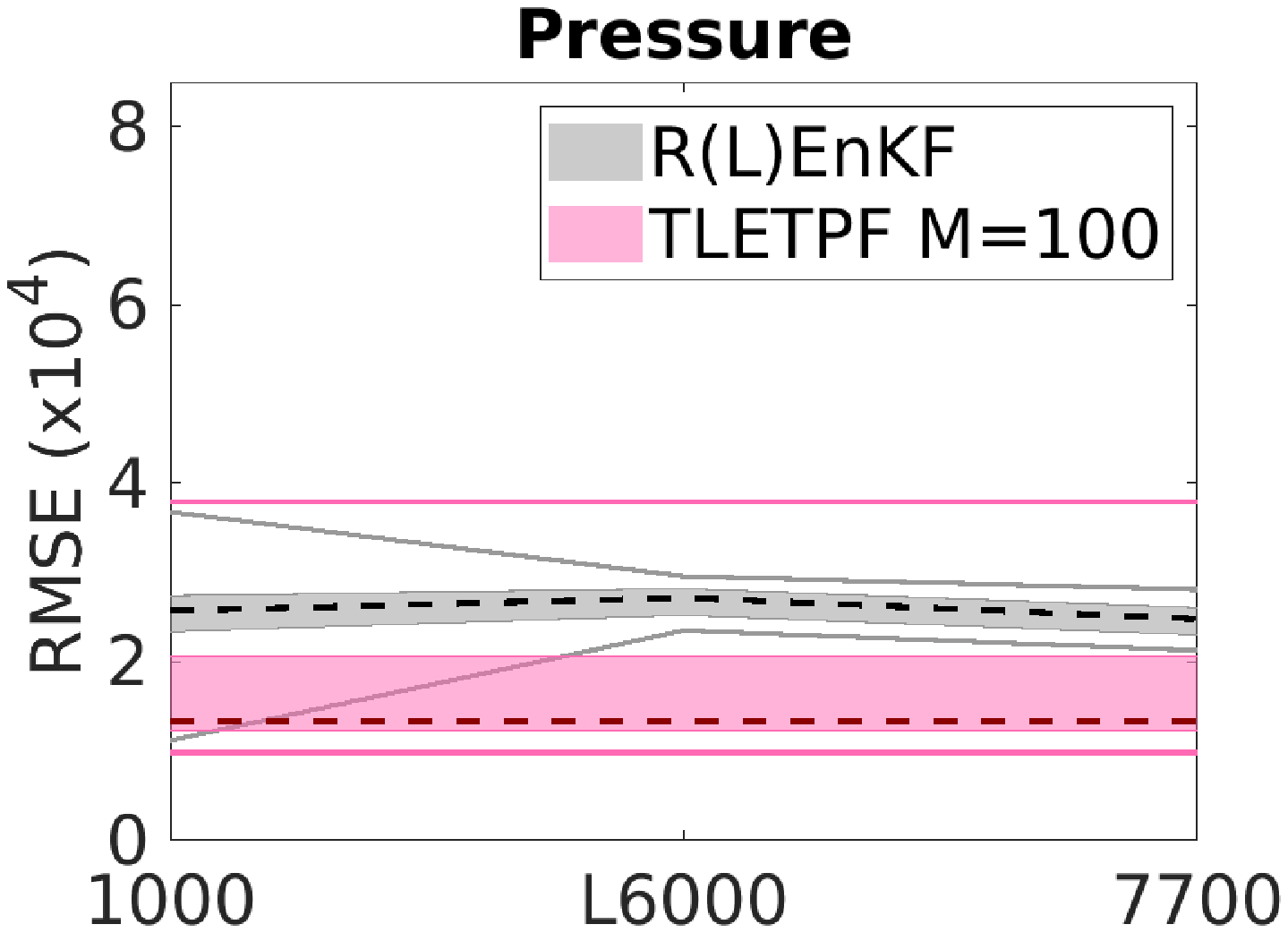}}
\caption{RMSE for $\log(\vec{k})$ is 
on the left. RMSE $\times10^4$ for pressure is on the right. 
A dashed line is for the median,
a shaded area is for 25 and 75 percentile, and a solid line is for 2 and 98 percentile over 20 simulations. 
On the x-axis numbers stand for ensemble sizes and L stands for a localized method.  
R(L)EnKF is shown in gray. TLETPF for ensemble size 100 is shown in pink.}\label{fig:Fig5}
\end{figure}

When comparing TLETPF to R(L)EnKF, we observe that the number of iterations is 2--3 times larger.
Let us fix computational cost for both methods in terms of model $g$ evaluations.
Then for an ensemble size $M=100$ of TLETPF, 
a computationally equivalent ensemble size of REnKF is 7700
and of RLEnKF is 6000.
Therefore we perform REnKF with ensemble size 7700
and RLEnKF with ensemble size 6000 (with favoring localization radius $r^{\rm loc}=6$).
We compute RMSE of mean log permeability and mean pressure and display them in Fig~\ref{fig:Fig5}
in gray for R(L)EnKF. TLETPF with ensemble size 100 is shown in pink in Fig~\ref{fig:Fig5}.
We observe that TLETPF at ensemble size 100 still gives better estimation of pressure 
than R(L)EnKF at an immense ensemble size. Therefore, 
the sensitivity to uncertain boundary conditions does not diminish
by extensively increasing ensemble size for R(L)EnKF. We should also note that 
though RLEnKF at ensemble size 1000 gives smaller errors on average than at ensemble size 6000,
errors vary more. This means RLEnKF at ensemble size 1000 is not robust. 

We have compared methods with respect to computational cost of evaluating the model $g$. 
Another computational cost is independent of the model evaluation and is associated with 
solving an optimization problem. For TLETPF it is $T^{\rm PF}(N^2M + M^2) \ln M$.
For RLEnKF it is $T^{\rm KF} M^2 \kappa$, assuming $M>\kappa$ and a low-rank approximation of covariance matrices~\cite{MaTo18}.
For TLETPF with $M=100$ computationally equivalent ensemble size of RLEnKF is around 500,
which is below 6000 considered above.
It is important to note, that $N^2$
in TLETPF can be distributed between computational nodes due to independence between the grid cells,
while in RLEnKF it cannot. Moreover, there also exists a computationally less expensive approximation, the Sinkhorn approximation~\cite{Acetal17}.


\section{Conclusions}
It has been known that EnKF is excellent in estimating Gaussian probabilities. 
Since log permeability is described by a Gaussian process, EnKF (iterative, regularized, or a smoother)
gives fine estimations of it. However, uncertainty in inverse problems is 
not only in rock properties but also in, for example, boundary conditions, geometry of the domain, and model simplifications. 
These sources of model error might not be described by Gaussian processes. 

Moreover, even if probability of model error is Gaussian or skewed but model error is non-additive,
it is not clear whether EnKF is able to give correct parameter-state estimation. 
In this paper, we have shown that R(L)EnKF is failing in such cases. 
For a low-dimensional test problem considered here, R(L)EnKF predicts well Gaussian distributed 
model error but gives poor estimation of multimode distributed parameter.
For a high-dimensional test problem of a Darcy flow with uncertain boundary conditions, 
R(L)EnKF predicts well Gaussian distributed parameter but gives poor estimation of skewed model error.
This results in inadequate estimation of pressure, which does not improve upon increasing ensemble size to 7700. 

T(L)ETPF, on the contrary, gives excellent estimation of both multimodal distribution for the low-dimensional test problem
and of skewed distribution for the high-dimensional test problem. 
For the high-dimensional test problem, even though log permeability estimation by T(L)ETPF is inferior to 
log permeability estimation by R(L)EnKF, pressure estimation is superior.
This is due to both well-estimated model error by T(L)ETPF and pressure sensitivity to uncertain boundary conditions.
Last but not least, T(L)ETPF with ensemble size 100 gives better pressure estimation 
than R(L)EnKF with ensemble size 7700 (6000).

\begin{acknowledgements}
We would like to acknowledge Marco Iglesias (U. of Nottingham) for fruitful discussions about 
sources of model error,
and for providing a code to model the Darcy flow and an MCMC code.
\end{acknowledgements}

\bibliographystyle{plain}   
\bibliography{DuRu18.bib}

\end{document}